\documentclass[12pt]{article}
\usepackage{amssymb}
\usepackage{amsmath,amsthm}
\usepackage[latin1]{inputenc}
\usepackage{graphicx}
\usepackage{epsfig}
\usepackage{hyperref}
\hypersetup{colorlinks=true, urlcolor=blue}

\newtheorem{Theorem}{\quad Theorem}[section]

\newtheorem{corollary}[Theorem]{\quad Corollary}

\newtheorem{proposition}[Theorem]{\quad Proposition}

\newtheorem{lemma}[Theorem]{\quad Lemma}

\newtheorem{remark}[Theorem]{\quad Remark}

\title{Offensive alliances in cubic graphs}

\author{\href{http://deim.urv.cat/~jarodriguez/}{J. A. Rodr\'{\i}guez}\footnote{e-mail:\mbox{\tt
    juanalberto.rodriguez\@@urv.net}} \\
{\em Department of Computer Engineering and Mathematics}\\
Rovira i Virgili University of Tarragona\\ Av. Pa\"{\i}sos Catalans
26, 43007 Tarragona, Spain \\
J. M. Sigarreta\footnote{e-mail:\mbox{\tt
    josemaria.sigarreta\@@uc3m.es}}\\
{\em Department of Mathematics }\\ Carlos III  University of Madrid\\
Avda. de la Universidad 30, 28911 Leganés (Madrid),  Spain }

\date{}

\begin{document}

\maketitle

\begin{abstract}
 An offensive alliance  in a graph $\Gamma=(V,E)$ is a set of
vertices $S\subset V$ where for every vertex $v$ in its boundary it
holds that the majority of vertices in $v$'s  closed neighborhood
are in $S$. In the case of strong offensive alliance, strict
majority is required. An alliance $S$ is called global if it affects
every vertex in $V\backslash S$, that is, $S$ is a dominating set of
$\Gamma$. The global offensive alliance number $\gamma_o(\Gamma)$
(respectively, global strong offensive alliance number
$\gamma_{\hat{o}}(\Gamma)$) is the minimum cardinality of a global
offensive (respectively, global strong offensive) alliance in
$\Gamma$. If $\Gamma$ has global independent offensive alliances,
then the \emph{global independent offensive alliance number}
$\gamma_i(\Gamma)$ is the minimum cardinality among all independent
global offensive alliances of $\Gamma$. In this paper we study
mathematical properties of the global (strong) alliance number of
cubic graphs. For instance, we show that for all connected cubic
graph of order $n$,
$$\frac{2n}{5}\le \gamma_i(\Gamma)\le  \frac{n}{2}\le  \gamma_{\hat{o}}(\Gamma)\le \frac{3n}{4}
\le \gamma_{\hat{o}}({\cal L}(\Gamma))=\gamma_{o}({\cal
L}(\Gamma))\le n,$$ where ${\cal L}(\Gamma)$ denotes the line graph
of $\Gamma$. All the above bounds are tight.
\end{abstract}



\section{Introduction}
The study of defensive alliances in graphs, together with a variety
of other kinds of alliances, was introduced in \cite{alliancesOne}.
In the cited paper there was initiated the study of the mathematical
properties of alliances. In particular, several bounds  on the
defensive alliance number were given. The particular case of global
(strong) defensive alliance was investigated  in
\cite{GlobalalliancesOne} where several bounds on the global
(strong) defensive alliance number were obtained.  In
\cite{spectral} there were obtained several tight bounds on
different types of alliance numbers of a graph, namely the (global)
defensive alliance number, (global) offensive alliance number and
(global) dual alliance number. In particular, there  was
investigated the relationship between the alliance numbers of a
graph and its algebraic connectivity, its spectral radius, and its
Laplacian spectral radius. A particular study of the alliance
numbers, for the case of planar  graphs, can be found in
\cite{planar}. For the study of offensive alliances we cite
\cite{chellali1,favaron,offensive} and for the study of alliances in
trees we cite \cite{chellali,chellali1,planar}. For the study of
alliance free sets and alliance cover sets we cite
\cite{kdaf,kdaf1,cover} and, finally, for the study of defensive
alliances in the line graph of a simple graph we cite \cite{line}.
The aim of this work is to study global offensive alliances in cubic
graphs.

\section{Preliminary Notes}

 In this paper
$\Gamma=(V,E)$ denotes a simple and connected graph of order $n$ and
${\cal L}(\Gamma)$ denotes the line graph of $\Gamma$. The degree of
a vertex $v\in V$ will be denoted by $\delta(v)$ and the subgraph
induced by a set $S\subset V$ will be denoted by $\langle S\rangle$.
For a non-empty subset $S\subset V$, and a vertex $v\in V$, we
denote by $N_S(v)$ the set of neighbors $v$ has in $S$:
$N_S(v):=\{u\in S: u\sim v\}$.
 Similarly, we denote by
$N_{V\setminus S}(v)$ the set of neighbors $v$ has in $V\setminus
S$: $N_{V\setminus S}(v):=\{u\in V\setminus S: u\sim v\}$.  The
boundary of a set $S\subset V$ is defined as $\partial
(S):=\displaystyle\bigcup_{v\in S}N_{V\setminus S}(v).$

A non-empty set of vertices $S\subset V$ is called {\em offensive
alliance} if and only if for every $v\in \partial (S)$, $| N_S(v) |
\ge | N_{V\setminus S}(v)|+1.$ That is, a non-empty set of vertices
$S\subset V$ is called  offensive alliance if and only if for every
$v\in \partial (S)$, $2| N_S(v) | \ge \delta(v)+1.$

An offensive alliance $S$ is called {\em strong} if for every vertex
$v\in \partial (S)$, $| N_S(v) | \ge | N_{V\setminus S}(v)|+2.$ In
other words, an offensive alliance $S$ is called  strong if for
every vertex $v\in \partial (S)$, $2| N_S(v) | \ge \delta(v)+2.$

The {\em   offensive alliance number} (respectively, \emph{strong
offensive alliance number}), denoted  $a_{o}(\Gamma)$ (respectively,
 $a_{\hat{o}}(\Gamma)$), is defined as
the minimum cardinality of an  offensive alliance (respectively,
strong offensive alliance) in $\Gamma$.

 A non-empty set of vertices $S\subset V$ is a {\em global offensive
alliance} if for every vertex $v\in V\setminus S$, $| N_S(v) | \ge |
N_{V\setminus S}(v)|+1$. Thus, global offensive alliances are also
dominating sets, and one can define the {\em global offensive
alliance number}, denoted $\gamma_{o}(\Gamma)$, to equal the minimum
cardinality of a global offensive alliance in $\Gamma$. Analogously,
$S\subset V$ is a {\em global strong offensive alliance} if for
every vertex $v\in V\setminus S$, $| N_S(v) | \ge | N_{V\setminus
S}(v)|+2,$ and the  {\em global strong offensive alliance number},
denoted $\gamma_{\hat{o}}(\Gamma)$, is defined as the minimum
cardinality of a global strong offensive alliance in $\Gamma$.

The \emph{independence  number} $\alpha(\Gamma)$ is the cardinality
of the largest independent set of $\Gamma$. A set is an
\emph{independent offensive alliance} in $\Gamma$ if it is an
offensive alliance and it is an independent set. If $\Gamma$ has
independent offensive alliances, then the \emph{independent
offensive alliance number} $a_i(\Gamma)$ is the minimum cardinality
among all independent offensive alliances of $\Gamma$. The
\emph{domination number} of $\Gamma$, denoted $\gamma(\Gamma)$, is
the minimum cardinality of a dominating set in $\Gamma$. The
\emph{independent domination number} $i(\Gamma)$ is the minimum
cardinality among all independent dominating sets of $\Gamma$.  If
$\Gamma$ has global independent offensive alliances, then the
\emph{global independent offensive alliance number}
$\gamma_i(\Gamma)$ is the minimum cardinality among all independent
global offensive alliances of $\Gamma$. Thus, $\gamma(\Gamma)\le
\gamma_o(\Gamma)\le \gamma_i(\Gamma)\le \alpha(\Gamma)$,
$a_i(\Gamma)\le \gamma_i(\Gamma)$ and $i(\Gamma)\le
\gamma_i(\Gamma)$.

In this paper we study mathematical properties of the global
(strong) alliance number of cubic graphs. For instance, we show that
for all connected cubic graph of order $n$,
$$\frac{2n}{5}\le \gamma_i(\Gamma)\le  \frac{n}{2}\le  \gamma_{\hat{o}}(\Gamma)\le \frac{3n}{4}
\le \gamma_{\hat{o}}({\cal L}(\Gamma))=\gamma_{o}({\cal
L}(\Gamma))\le n.$$  All the above bounds are tight.

\section{Main Results} These are the main results of the paper.

\begin{Theorem} \label{th1}
Let $\Gamma$ be a connected cubic graph of order $n$.
\begin{enumerate}
\item $\frac{n}{2}\le \gamma_{\hat{o}}(\Gamma)\le \frac{3n}{4}$.
\item $\gamma_{\hat{o}}(\Gamma)= \frac{n}{2}$ if and only if $\Gamma$
is a bipartite graph.

\item $\gamma_{\hat{o}}(\Gamma)= \frac{3n}{4}$ if and only if $\Gamma$
is isomorphic to the complete graph $K_4$.
\end{enumerate}
\end{Theorem}

\begin{proof}

\mbox{  }

\begin{enumerate}
\item If $\Gamma=(V,E)$ is $\delta$-regular, then $S\subset V$ is a
$\delta$-dominating set if and only if $V\backslash S$ is an
independent set. Therefore,
\begin{equation}
\gamma_\delta(\Gamma)+\alpha(\Gamma)=n.
\end{equation}
 Moreover, if $\Gamma$
is a cubic graph,  a set $S\subset V$ is a strong offensive alliance
if and only if $S$ is a 3-dominating set. Therefore,
\begin{equation}\label{strongIndep}
\gamma_{_3}(\Gamma)=\gamma_{\hat{o}}(\Gamma)=n-\alpha(\Gamma).
\end{equation}
Finally, as for all $\delta$-regular graph $\Gamma$,
$\alpha(\Gamma)\le \frac{n}{2}$, the bound follows.

On the other hand, for all  global strong offensive alliance $S$
such that $|S|=\gamma_{\hat{0}}(\Gamma)$, $V\backslash S$ is an
independent set. Thus, $\frac{3n}{2}\leq
3(n-\gamma_{\hat{0}}(\Gamma))+\gamma_{\hat{0}}(\Gamma)$. Hence, the
upper bound follows.

\item
If $\Gamma=(V,E)$ is a bipartite cubic graph, then each set of the
bipartition of $V$ is a strong global offensive alliance in $\Gamma$
of cardinality $\frac{n}{2}$, so
$\gamma_{\hat{0}}(\Gamma))=\frac{n}{2}$. Conversely, if $S$ is a
$3$-dominating set of cardinality $\frac{n}{2}$ in a cubic graph,
then the edge-cut between $S$ and $V\backslash S$ has size
$3\frac{n}{2}$. Hence, both, $S$ and $V\backslash S$, are
independent sets, so $\Gamma$ is a bipartite graph.

\item We only need to show that  $\alpha(\Gamma)= \frac{n}{4} \Rightarrow \Gamma \cong K_4$.
 Suppose $\Gamma\ncong K_4$. If $X$ is an independent set such
that $|X|=\alpha(\Gamma)=\frac{n}{4}$, then $\langle V\backslash
X\rangle$ is the disjoin union of cycles. Let $x_i$ be the number of
cycles of length $i$ in $ V\backslash X$, $i=3,...,\frac{3n}{4}$.
Thus,
$$\sum ix_i=\frac{3n}{4}$$
If $x_i>0$, for some $i>3$, then
$$\sum 3\left\lfloor\frac{i}{2}\right\rfloor x_i>\sum i x_i=\frac{3n}{4}.$$
Thus,
$$\alpha\left( \langle V\backslash X\rangle \right)=\sum \left\lfloor\frac{i}{2}\right\rfloor x_i>\frac{n}{4}.$$
Therefore, $x_i=0$, for all $i>3$. Let $Y_i\subset V\backslash X$
such that $\langle Y_i\rangle\cong K_3$, $i=3,...,\frac{n}{4}$. Let
$x\in X$. As $\Gamma$ is connected and $\Gamma \ncong K_4$, $\langle
Y_i \cup \{ x \}\rangle\ncong K_4$ and, in consequence, $\alpha
\left(\langle Y_i \cup \{ x \}\rangle\right)=2$, $\forall i$. Let
$W=\{x,y_1,y_2,...,y_{\frac{n}{4}}\}$, where $y_i\in Y_i$ and
$y_i\nsim x$. As $W$ is an independent set, the result follows by
contradiction.
\end{enumerate}
\end{proof}

It is well-known that the independent set problem is NP-complete
\cite{Garey}. Hence, a direct consequence of (\ref{strongIndep}) is
the following:

\begin{remark}The minimum global strong offensive alliance problem is NP-complete.
\end{remark}

 There are some classes of cubic graphs in which we can compute the global strong offensive alliance number in terms
 of the order. For instance, if $\Gamma={\cal L}(\Gamma_1)$ is a cubic graph of order $n\ge 6$, then we
 have
$\gamma_{\hat{0}}(\Gamma)=2\gamma_{\hat{0}}(\Gamma_1)=\frac{2n}{3}.$
That is, $\Gamma_1=(V,E)$ is a bipartite semiregular  graph of
degrees $2$ and $3$. If  $V_{3}\subset V$ denotes the vertex-set of
degree 3 in $\Gamma_1$, then
$\gamma_{\hat{0}}(\Gamma_1)=|V_{3}|=\alpha(\Gamma)$. On the other
hand, $|V_{3}|=\frac{n}{3}$. Therefore,
$\gamma_{\hat{0}}(\Gamma)=n-\alpha \left(\Gamma)\right)=
3\gamma_{\hat{0}}(\Gamma)-\gamma_{\hat{0}}(\Gamma)=\frac{2n}{3}$.
Notice that in this class of graphs are included the cubic graphs of
order $n\ge 6$ having $\frac{n}{3}$ independent triangles\footnote{A
set of  triangles    is independent if it contains no common
vertices.}.

\begin{Theorem}\label{thOf}  {\rm \cite{spectral,offensive}}
Let $\Gamma$ be a simple graph of order $n$ and minimum degree
$\delta$. Let $\mu$ be the Laplacian spectral radius of $\Gamma$.
Then
$$\frac{n}{\mu}\left\lceil\frac{\delta +1}{2}\right\rceil \le
\gamma_{o}(\Gamma) \le \displaystyle\frac{n(2\mu-\delta)}{2\mu}.$$
\end{Theorem}

\begin{corollary} \label{coro}
Let $\Gamma$ be a $\delta$-regular graph of order $n$. Then
$$\frac{n}{4}\left\lceil \frac{2\delta-1}{2}\right\rceil \le \gamma_{o}({\cal L}(\Gamma))\le \frac{n(\delta+1)}{4}.$$
\end{corollary}

\begin{proof}
We denote by $A$ the adjacency matrix of ${\cal L}(\Gamma)$ and by
$2(\delta-1)=\lambda_0>\lambda_1>\cdots>\lambda_b=-2$ its distinct
eigenvalues. We denote by $L$ the Laplacian matrix of ${\cal
L}(\Gamma)$ and by $\mu_0=0<\mu_1<\cdots <\mu_b$ its distinct
Laplacian eigenvalues. Then, since $L= 2(\delta-1)I_n- A$, the
eigenvalues of both matrices, $A$ and $L$, are related by
\begin{equation} \label{eigen}
\mu_l = 2(\delta-1)-\lambda_l, \quad l=0, \dots, b.
\end{equation}
Thus, the Laplacian spectral radius of ${\cal L}(\Gamma)$ is
$\mu_b=2\delta$. Therefore, the result immediately follows
\end{proof}

In the case of regular graphs of even degree all global offensive
alliance is strong. Hence we have
$\gamma_{o}(\Gamma)=\gamma_{\hat{o}}(\Gamma)$. Therefore, by
Corollary \ref{coro} we obtain that if $\Gamma$ is a cubic graph,
then
\begin{equation}
\frac{3n}{4}\le \gamma_{\hat{o}}({\cal L}(\Gamma))=\gamma_{o}({\cal
L}(\Gamma))\le n.
\end{equation}

The above bounds are tight. Example of equality in above upper bound
is the complete graph of order 4: $\gamma_{o}({\cal L}(K_4)=4=n$. In
the case of the bipartite complete graph $\Gamma=K_{3,3}$  we have $
\gamma_{o}({\cal L}(K_{3,3}))=5$ and the above lower bound gives
$\frac{9}{2}\le \gamma_{o}({\cal L}(K_{3,3}))$.

\begin{Theorem}{\rm \cite{offensive}} \label{cotainf}
For all  graph $\Gamma$ of order $n$, minimum degree $\delta$ and
maximum degree $\Delta$,
 $$
 \frac{2n}{3} \ge\gamma_0(\Gamma)\ge \left\lbrace \begin{array}{ll}
\left\lceil\frac{n(\delta+1)}{2\Delta+\delta+1}\right\rceil  & {\rm
if }\quad \delta \quad {\rm odd};
                            \\
                            \\
                            \left\lceil\frac{n\delta}{2\Delta+\delta}\right\rceil  &  {\rm otherwise;} \end{array}
                                                \right  .$$

\end{Theorem}

 As we show in the following
result, the bound $\gamma_o(\Gamma)\le  \frac{2n}{3}$ is improved
for the case of regular graphs of odd degree.

\begin{Theorem}  \label{cotasup-n/2}
For all regular graph $\Gamma$  of order $n$ and odd degree
$\delta$, $$\frac{n(\delta+1)}{3\delta+1}\le \gamma_o(\Gamma)\le
\frac{n}{2}.$$
\end{Theorem}

\begin{proof}
The lower bound is a particular case of Theorem \ref{cotainf}. Let
$\{X, Y\}$ be a bipartition of $V$ such that $|X|=|Y|=\frac{n}{2}$
and the edge-cut between $X$ and $Y$ has maximum size. If neither
$X$ nor  $Y$ are global offensive alliances in $\Gamma$, then there
exist $x\in X$ and $y\in Y$ such that $|N_Y(x)|\le |N_X(x)|$ and
$|N_X(y)|\le |N_Y(y)|$. Since $\Gamma$ is regular of odd degree,
$|N_Y(x)| < |N_X(x)|$ and $|N_X(y)| < |N_Y(y)|$. Thus, by
contradiction we deduce the result.
\end{proof}

In the case of cubic graphs we have,
\begin{equation}
\frac{2n}{5}\le \gamma_o(\Gamma)\le  \frac{n}{2}.
\end{equation}

Example of equality in above upper bound is the family of cubic
graphs $\Gamma=C_r\times K_2$, where $C_r$ denotes the $r$-cycle
graph. In this case $\gamma_o(\Gamma)=r$.

\begin{proposition}
Let $\Gamma$ be a cubic graph of order $n$. All global offensive
alliance in $\Gamma$ of cardinality $\frac{2n}{5}$ is an independent
set.
\end{proposition}

\begin{proof}
Let $X\subset V$ be a global offensive alliance in $\Gamma$ of
cardinality $\frac{2n}{5}$. Let $c$ be the size of the edge-cut
between $X$ and $V \backslash X$.As $X$ is a $2$-dominating set,
$2\frac{3n}{5} \le c$. Moreover, $c\le 3\frac{2n}{5} $, so the size
of $\langle X \rangle$ is zero.
\end{proof}

A set nonempty set $S\subseteq V$ is a strong defensive alliance in
$\Gamma$  if for every $v\in S$, $2|N_S(v)|\ge \delta(v)$. A set
$X\subseteq V$ is  {\em strong defensive alliance free} if for all
strong defensive alliance $S$, $S\setminus X\neq\emptyset$, i.e.,
$X$ do not contain any strong defensive alliance as a subset
\cite{kdaf}. A strong defensive alliance free set $X$ is
\emph{maximal} if  for all $v\notin X$, exists $S\subseteq X$ such
that $S\cup \{v\}$ is a strong defensive alliance. A \emph{maximum}
strong defensive alliance free set is a maximal  strong defensive
alliance free set of largest cardinality. We denote by
$\phi_0(\Gamma)$   the cardinality of a maximum strong defensive
alliance free set of $\Gamma$.

Similarly, a set $Y \subseteq V$ is a \emph{strong defensive
alliance cover} if for all strong defensive alliance $S$, $S\cap
Y\neq\emptyset$, i.e., $Y$ contains at least one vertex from each
strong defensive alliance of $\Gamma$.  A strong defensive alliance
cover $Y$ is \emph{minimal} if no proper subset of $Y$ is a strong
defensive alliance cover. A \emph{minimum}  strong defensive
alliance cover is a minimal cover of smallest cardinality. We will
denote by $\zeta_0(\Gamma)$
  the cardinality of a minimum  strong
defensive alliance cover of $\Gamma$.

\begin{lemma}{\rm \cite{cover}}  \label{lemma1}
If  $X\subset V$ is a global offensive alliance in $\Gamma=(V,E)$,
then the set $V\setminus X$ is strong defensive alliance free in
$\Gamma$.
\end{lemma}

\begin{lemma}{\rm \cite{kdaf}} \label{lemadutton}
If each block\footnote{A block is a maximal 2-connected subgraph of
a given graph $\Gamma$.} of $\Gamma$ is an odd clique or an odd
cycle, then $\phi_0(\Gamma)\le \zeta_o(\Gamma)$.
\end{lemma}

\begin{Theorem}
Let $\Gamma$ be a cubic graph of order $n$. If each block of
$\Gamma$ is  an odd cycle, then $$ \gamma_o(\Gamma)=\frac{n}{2}.$$
\end{Theorem}

\begin{proof}
By Lemma \ref{lemadutton} and $\phi_{0}(\Gamma)+\zeta_0(\Gamma)=n$
we have $\phi_0(\Gamma)\le \frac{n}{2}$. Moreover, by Lemma
\ref{lemma1} we have  $n-\gamma_o(\Gamma)\le \phi_0(\Gamma)$. Hence,
$\gamma_o(\Gamma)\ge \frac{n}{2}$. By Theorem \ref{cotasup-n/2} we
conclude the proof.
\end{proof}

\subsection{Independent offensive alliances}\label{sectionIndep}

As a consequence of Theorem \ref{th1} we obtain the following
result.

\begin{corollary}
 Let $\Gamma$ be a cubic graph. $\Gamma$ has a global strong independent offensive
alliance if and only if $\Gamma$ is a bipartite graph.
\end{corollary}

 \begin{proof}In the case of cubic graphs we have $\frac{n}{2}\le
\gamma_{\hat{o}}(\Gamma)$ and $\alpha(\Gamma)\le\frac{n}{2}$. Both
equalities holds true if and only if $\Gamma$ is a bipartite graph.
\end{proof}

In the case of a bipartite cubic graph we have,
$$\gamma_i(\Gamma)=\gamma_{o}(\Gamma)=\gamma_{\hat{o}}(\Gamma)=\gamma(\Gamma)=\alpha(\Gamma)=\frac{n}{2}.$$

For all graph having independent offensive alliances we have $
\gamma_o(\Gamma)\le  \gamma_i(\Gamma)$. Example of equality is the
Petersen graph $O_3$ and the class of bipartite cubic graphs.
Obviously, if $\alpha(\Gamma)< \gamma_o(\Gamma)$, then $\Gamma$ do
not contains global independent offensive alliances. Examples of
graphs having $\alpha(\Gamma)<\gamma_o(\Gamma)$ are the
 graphs isomorphic to $C_{2k+1}\times K_2$. In this case
 $2k=\alpha(\Gamma)<\gamma_o(\Gamma)=2k+1$.

\begin{Theorem}
Let $\Gamma=(V,E)$ be a cubic graph of order $n$. If $\Gamma$ has a
global independent offensive alliance, then:
\begin{enumerate}
\item $\frac{2n}{5}\le \gamma_i(\Gamma)\le  \frac{n}{2}.$

\item $\gamma_i(\Gamma)=\frac{n}{2}$ if and only if $\Gamma$ is
a bipartite graph.

\item  $\gamma_i(\Gamma)=\frac{2n}{5}$ if and only if there exists an
independent set $X\subset V$ such that $\langle V\backslash X
\rangle$ is a 1-factor of size $\frac{3n}{10}$.
\end{enumerate}
\end{Theorem}

\begin{proof}
Let $\Gamma=(V,E)$ be a cubic graph of order $n$.
\begin{enumerate}
\item If $\Gamma$ has a
global independent offensive alliance, then $\frac{2n}{5}\le
\gamma_o(\Gamma)\le \gamma_i(\Gamma)\le \alpha(\Gamma) \le
\frac{n}{2}.$

\item If  $X\subset V$ is an independent set, then
 $V\backslash X$ is  a  3-dominating set of $\Gamma$. Hence, if  $|X|=\frac{n}{2}$, the edge-cut
  between $X$ and $V\backslash X$ has size $3\frac{n}{2}$,
 so $\Gamma$ is a bipartite graph. Conversely, if $\Gamma$ is a
 bipartite cubic
 graph, then $\gamma_i(\Gamma)=\frac{n}{2}$.

\item  Let  $X$ be a global independent offensive alliance in $\Gamma$ such that
$|X|=\frac{2n}{5}$. The edge-cut
  between $X$ and $V\backslash X$ has size $3\frac{2n}{5}$. Since
$|V\backslash X|=\frac{3n}{5}$ and $X$ is a $2$-dominating set, each
vertex of
 $V\backslash X$ have one neighbor in $V\backslash X$,  so $\langle V\backslash X \rangle$
 is a 1-factor of size $\frac{3n}{10}$. The converse is immediate.
\end{enumerate}
\end{proof}

\begin{corollary}
Let $\Gamma$ be a cubic graph of order $n$. If
$\alpha(\Gamma)<\frac{2n}{5}$, then $\Gamma$ do not contains global
independent offensive alliances.
\end{corollary}

Examples of graphs having $\alpha(\Gamma)<\frac{2n}{5}$ are the
complete graph, $K_4$, and the graph $\Gamma=K_{3}\times K_2$.

\begin{Theorem} \label{otfen22}
Let $\Gamma$ be  a cubic graph of order $n$.  If $a_i(\Gamma)<
\gamma_i(\Gamma)$, then
$$\frac{n+2}{4}\leq a_i(\Gamma)\leq \frac{n-2}{2}.$$
\end{Theorem}

\begin{proof}
If   $S\subset V$ is an independent set in $\Gamma=(V,E)$, then
\begin{equation}\label{ofensiva}
\frac{3n}{2}=3|S|+\frac{\sum_{v\in V\backslash S}|N_{V\backslash
S}(v)|}{2}.
\end{equation}
Moreover, if $S\subset V$ is an independent offensive alliance in
$\Gamma$ such that $|S|=a_i(\Gamma)<\gamma_i(\Gamma)$, then there
exists at least one vertex $v\in V$ such that $v\notin (S\cup
\partial S)$.
Hence,
\begin{equation}\label{ofena}
\frac{\sum_{v\in V\backslash S}|N_{V\backslash S}(v)|}{2}\geq 3.
\end{equation}
By (\ref{ofensiva}) and (\ref{ofena}) we obtain the upper bound. On
the other hand, since $S$ is an offensive alliance in
$\Gamma=(V,E)$, then
\begin{equation}\label{ofensiva1}
3|S|\geq\sum_{v\in \partial S}|N_{S}(v)|\geq\sum_{v\in \partial
S}|N_{V\backslash S}(v)|+|\partial S|.
\end{equation}
Moreover, as $S$ is an independent set $|\partial S|\geq 3.$ Thus,
\begin{equation}\label{ofensiva2}
3|S|+\sum_{v\in
\partial S}|N_{V\backslash S}(v)|\ge \frac{3n}{2} .
\end{equation}
By (\ref{ofensiva}) and (\ref{ofensiva1}) we obtain de result.

\end{proof}

Example of equality in above bounds is the complete bipartite graph
$\Gamma=K_{3,3}$. In this case $a_i(\Gamma)=2$.


\end{document}